\documentclass[12pt,reqno]{amsart}

\usepackage[dvips]{graphicx} 

\usepackage
{hyperref}

\usepackage{breakurl}   

\theoremstyle{definition}

\numberwithin{equation}{section}

\newcommand\N {{\mathbb N}} 

\newcommand\R {{\mathbb R}}

\newcommand\astr {{}^\ast{\mathbb R}}

\newcommand\astn {{}^\ast{\mathbb N}}

\title {A Leibniz/NSA comparison}

\author[M. Katz]{Mikhail G. Katz} \address{Department of Mathematics,
  Bar Ilan University, Ramat Gan 5290002 Israel}
\email{katzmik@math.biu.ac.il}

\author[K. Kuhlemann]{Karl Kuhlemann}\address{Gottfried Wilhelm
  Leibniz University Hannover, D-30167 Hannover, Germany}
\email{kus.kuhlemann@t-online.de}

\author[D. Sherry]{David Sherry} \address{Department of Philosophy,
  Northern Arizona University, Flagstaff, AZ 86011, US}
\email{David.Sherry@nau.edu}

\begin{document}

\begin{abstract}
We present some similarities between Leibnizian and Robinsonian
calculi, and address some objections raised by historians.  The
comparison with NSA facilitates our appreciation of some Leibnizian
procedures that may otherwise seem obscure.  We argue that Leibniz
used genuine infinitesimals and infinite quantities which are not
merely stenography for Archimedean Exhaustion and that Leibniz's
procedures therefore find better proxies in NSA than in modern
Weierstrassian mathematics.
\end{abstract}

\thispagestyle{empty}


\maketitle


\section{Leibniz and NSA}
\label{s1}

The Leibnizian infinitesimal calculus represents a different approach
to the calculus from the one generally adopted today.  A seminal study
of Leibnizian methodology was published by Bos \cite{Bo74}.  A number
of distinct interpretations of the Leibnizian calculus have been
pursued by Leibniz scholars over the past few decades.

The article \cite{22b}, which analyzed a pair of rival interpretations
of the Leibnizian calculus, mentioned neither Robinson nor Nonstandard
Analysis/NSA (though the article made a passing mention of nonstandard
models of arithmetic, which go back to Skolem in 1933; see
\cite{Sk33}).  The present text seeks to make the Leibniz/NSA
comparison more explicit.  We mainly follow the interpretation of
Leibniz given by Bos~\cite{Bo74}.

In Section~\ref{two}, we summarize some relevant procedures of the
Leibnizian calculus, and signal some similarities with the Robinsonian
calculus.  The similarities are summarized in Table~\ref{tt}.  In
Section~\ref{s2}, we deal with some objections by historians.

\noindent
\begin{table}
\renewcommand{\arraystretch}{1.5}
\begin{tabular}[t]  
{|
@{\hspace{3pt}}p{2.4in}|| 
@{\hspace{3pt}}p{2.4in}|
} 
\hline Leibniz & NSA  
\\ 
\hline\hline
relation of infinite proximity & standard part (shadow)  
\\
\hline 
law of continuity & transfer principle  
\\
\hline 
assignable vs inassignable number & standard vs nonstandard number
\\
\hline 
infinitum terminatum & line segment of unlimited length
\\
\hline 
infinitesimal violates Euclid V.4 & infinitesimal violates the
Archimedean property
\\
\hline 
constant differentials & uniform partition into infinitesimal
subsegments
\\ 
\hline
\end{tabular}
\renewcommand{\arraystretch}{1}
\caption{\textsf{Leibniz and NSA${}^{\phantom{I}}$}}
\label{tt}
\end{table}

\section{Six similarities}
\label{two}

We present six similarities between the Leibnizian and Robinsonian
calculi, as summarized in Table~\ref{tt}.

\subsection{What kind of equality?}

Leibniz mentioned that he worked, not with the relation of
\emph{exact} equality between quantities, but rather equality up to
negligible infinitesimal terms, as when he wrote:
\begin{enumerate}\item[]
I think that not only those things are equal whose difference is
absolutely zero, but also those whose difference is incomparably
small.  And although this [difference] need not absolutely be called
Nothing, neither is it a quantity comparable to those whose difference
it is.%
\footnote{Leibniz \cite[p.\;322]{Le95b}.}
\end{enumerate}
This can be described as a relation of \emph{infinite proximity}.  In
NSA, a related procedure involves extracting the \emph{standard part}
of a limited quantity.%
\footnote{A number is \emph{limited} if it is smaller in absolute
  value than some standard number.}

\subsection{Law of Continuity}

One of the formulations of Leibniz's \emph{law of continuity} is ``the
rules of the finite are found to succeed in the infinite and vice
versa.''%
\footnote{``il se trouve que les regles du fini reussissent dans
  l'infini, {\ldots} et que vice versa les regles de l'infini
  reussissent dans le fini'' Leibniz \cite[pp.\;93--94] {Le02}.
  Cf.\;Robinson \cite[p.\;266]{Ro66}.}
In reference to Leibniz's law, Robinson asserted the following:
\begin{enumerate}\item[]
this is remarkably close to our transfer of statements from~$\R$
to~${}^\ast\R$ and in the opposite direction.%
\footnote{Robinson \cite[p.\;266]{Ro66}.  The transfer principle is a
  theorem asserting that a formula that holds for all standard inputs
  would in fact hold for all inputs.  Thus, if we had a proof that a
  relation such as~$\cos^2 x +\sin^2 x=1$ holds for all standard~$x$,
  it would automatically hold for all~$x$ by the transfer principle.}
\end{enumerate}

\subsection{Assignable or inassignable?}

Leibniz mentioned frequently a distinction between assignable and
inassignable numbers, as when he writes: 
\begin{enumerate}\item[]
Here~$dx$ means the element, that is, the (instantaneous) increment or
decrement, of the (continually) increasing quantity~$x$.  It is also
called difference, namely the difference between two proximate~$x's$
which differ by an element (or by an \emph{inassignable}), the one
originating from the other, as the other increases or decreases
(momentaneously).%
\footnote{Leibniz as translated by Bos \cite[p.\,18]{Bo74}; emphasis
  added.}
\end{enumerate}
In NSA, a parallel distinction exists between standard and nonstandard
numbers.  The distinction is used routinely in the procedures of
modern infinitesimal analysis.  Thus, an infinitesimal is defined as a
number smaller in absolute value than every standard positive real
number, and the continuity of a function~$f$ at a standard point~$x$
is defined by the condition that an infinitesimal change~$\alpha$ of
the independent variable always produces an infinitesimal change
$f(x+\alpha)-f(x)$ of the function.

\subsection{What is \emph{bounded infinity}?}
\label{s24}

Starting in his 1676 \emph{De Quadratura Arithmetica} \cite{Le16},
Leibniz used a type of line segment he referred to as \emph{infinitum
  terminatum}, which literally means ``bounded infinity''.  This is
greater than any assignable line segment.  What could the expression
\emph{bounded infinity} possibly mean?  In his chapter on the history
of mathematics \cite[Chapter 10]{Ro66}, Robinson makes no mention of
Leibniz's use of such an expression.  To a modern reader trained in
the Weierstrassian paradigm, the expression \emph{bounded infinity}
may seem like a contradiction in terms.  The pioneering work of
Cantor, Weierstrass, and others in creating rigorous foundations for
analysis has conditioned successive generations of mathematicians to
hold that there are just two types of infinite number in mathematics:
\emph{ordinals} and \emph{cardinals}.%
\footnote{To Leibniz, such concepts of modern set theory would
  arguably have been meaningless, since he rejected the notion of an
  \emph{infinite whole} as contradictory.  To Leibniz, \emph{bounded
    infinity} was the consistent counterpart of \emph{unbounded
    infinity}, such as the unbounded real line in modern terminology.
  For further analysis, see \cite{21a}.}
It is helpful to go beyond that dichotomy, and allow also for the
possibility of what could be called \emph{ringinals}; namely, elements
in a number system that are greater than all naive integers (such as a
nonstandard integer in the semiring~$\astn$ extending~$\N$).  Such
numbers are referred to as \emph{unlimited}.

Leibniz's notion of \emph{infinitum terminatum} is analogous to
segments of unlimited length (on a nonstandard line), greater than
segments of any standard length.  The notion played a crucial role in
the Leibnizian calculus, by enabling non-Archimedean phenomena (see
Section~\ref{s15}).  Thus, an infinitesimal is 
\begin{enumerate}\item[]
[an] infinitely small
fraction, or one whose denominator is an infinite number.%
\footnote{``fraction infiniment petite, ou dont le denominateur soit
  un nombre infini'' (Leibniz \cite[p.\;93]{Le02}).}
\end{enumerate}
The rules for operating with infinitesimal differentials, and their
applications to developing the basic procedures of the calculus,
including the Leibniz rule for the product, first appeared in
Leibniz's 1684 article~\cite{Le84}.

\subsection{Leibniz on violation of Euclid V, Definition 4}
\label{s15}

In a 1695 article \cite{Le95b} containing a rebuttal of a criticism by
Nieuwentijt, Leibniz made it clear that his ``incomparable''
quantities, when compared to 1, violate the comparability notion
formulated in Euclid Book~V, Definition\;5 (appearing as Definition\;4
in modern editions). Euclid's definition reads as follows:
\begin{quote}
Any magnitude may be multiplied as many times as to exceed any given
homogeneous magnitude.%
\footnote{De Risi \cite[p.\;626]{De16}.}
\end{quote}
Leibniz made similar comments in a famous 1695 letter to l'Hospital.%
\footnote{Leibniz \cite[p.\;288]{Le95a}.}
Similarly, NSA infinitesimals violate the Archimedean property.

\subsection{Constant differentials}

Leibniz historians often mention his use of terms such as
``progression of variables'' and ``constant differentials.''  Thus,
Bos writes: 
\begin{enumerate}\item[]
Leibniz explains that if~$dx$ is taken constant, one may treat the
quadrature as~$\int y$ (`sum of all~$y$'), as is done in the theory of
indivisibles, but if one wishes to consider different progression of
the variables, the quadrature has to be evaluated as~$\int y\,dx$.%
\footnote{Bos \cite[p.\;27]{Bo74}.}
\end{enumerate}
Such comments can be formalized in modern infinitesimal analysis in
terms of a uniform partition of the segment into infinitesimal
subsegments, when the integral can be computed via the Riemann sum
corresponding to the uniform partition.%
\footnote{Of course, the distinction between uniform and non-uniform
  partitions of the interval of integration can equally well be
  formulated in traditional non-infinitesimal analysis; however, the
  latter lacks the language to express partitions into
  \emph{infinitesimal} subsegments, which would necessarily have to be
  paraphrased in terms of limits of sequences of Archimedean
  partitions.}
Then ``nonconstant differentials'' would correspond to partitions that
are not necessarily uniform.  The idea of ``progression of the
variables'' may seem odd to a reader accustomed to thinking of a
variable as ranging through a continuous range such a subinterval
of~$\R$, but it is natural from the viewpoint of infinitesimal
partitions exploited in NSA.

\section{Possible objections}
\label{s2}

Among historians of mathematics, critical attitudes are not uncommon
toward attempts to mention Leibniz and NSA in the same article.  In
the context of a discussion of the Leibnizian calculus, a noted
historian expressed the following reservation in march 2023 in a
private communication: ``the view that NSA is pervasive in the history
of math is a legitimate (though wrong, to my view) thesis that anyone
has the right to defend.''%
\footnote{Marco Panza, private communication, 17 march 2023.}
There are at least three possible ways of interpreting the claim that
would be the intended target of such a reservation:
\begin{enumerate}
\item
\label{c1}
Leibniz used ultrafilters;
\item
\label{c2}
Leibniz's procedures find better proxies in NSA than in modern
Weierstrassian mathematics;
\item
\label{c3}
Leibniz's infinitesimals were non-Archimedean, and not merely
stenography for Archimedean exhaustion.
\end{enumerate}

Clearly, claim \eqref{c1} would be an anachronism, whether it is made
with regard to ultrafilters and ultrapowers or other modern
set-theoretic constructions of nonstandard number fields.%
%
%

With regard to claim \eqref{c2}, it would perhaps be reasonable to
maintain that Leibniz was neither Weierstrassian nor Robinsonian,
while both can be helpful for understanding aspects of his work.
Nevertheless, one can analyze the \emph{procedures} of the Leibnizian
calculus to gauge which modern theory provides better proxies for the
kind of inferential moves used by Leibniz.  Thus, the importance of
the notion of the \emph{infinitum terminatum} in the Leibnizian
calculus and geometry (see Section~\ref{s24}) was emphasized in the
early articles by Leibniz historian Knobloch \cite{Kn90}, \cite{Kn94}.
As mentioned in Section~\ref{s1}, there is a ready analog in
nonstandard analysis, given by segments of unlimited length such as a
nonstandard integer.  Finding proxies in a Weierstrassian framework
would be more difficult, and necessarily involve an infinite process
rather than a single segment.

With regard to claim~\eqref{c3}, certainly Leibniz was no mere
stenographer.  More to the point, one can ask whether, whenever
Leibniz used the term \emph{infinitesimal}, he just meant the reader
to fill in the details of the corresponding Archimedean exhaustion
argument, as in an Ishiguroan interpretation%
\footnote{Ishiguro \cite[Chapter~5]{Is90}.}
which reads the passages containing the term \emph{infinitesimal} as
shorthand for Archimedean exhaustion; this particular debate was
already analyzed in \emph{Depictions} \cite{22b}.  Leibniz scholarship
today is an area of lively debate; see e.g., Archibald et
al.~\cite{Ar22}, Bair et al.~\cite{22a}, \cite{23a}, Ugaglia and Katz
\cite{Ug23}, Katz et al.~\cite{23h}, \cite{24c}.

We comment briefly on possible interpretations of the work of other
historical mathematicians in NSA terms.  While Cauchy often used
genuine infinitesimals, as detailed in \cite{20b} and \cite{22a}, he
had less interest in philosophy than Leibniz.  Cauchy produced
mathematics for the most part without looking for underlying
philosophical principles (with the possible exception of his Torino
lectures \cite{Ca33}).  While there are a few mentions of assignable
quantities in Cauchy, one does not find any mention of inassignables,
and certainly no occurrences of \emph{infinita terminata}.  Leibniz
viewed principles such as law of continuity as philosophical insights.
Today we view its counterpart, the transfer principle, as a
mathematical insight.  Cauchy wasn't interested in such general
principles as much as Leibniz, and therefore, from a modern viewpoint,
did not get as far as Leibniz foundationally speaking, even though he
was of course much ahead of Leibniz mathematically speaking.  As far
as Fermat is concerned, he was extremely cautious and we have little
indication of what he thought about his increments~$E$.  Indeed, there
were good reasons for his caution; see \cite{18d}.  We have argued
that adequality is best understood as some kind of relation of
infinite proximity; see also \cite{18d}.

\section{Conclusion}

The viability of interpreting the Leibnizian calculus in terms of NSA
depends crucially on the distinction between procedures and
foundations.  The analogies we developed here concern the realm of
procedures (such as those listed in Table~\ref{tt}).  As far as
foundations are concerned, it may be pertinent to mention the
existence of (at least) two distinct approaches to NSA.  Robinson's
approach as developed in the 1960s involved extending the field of
real numbers~$\R$ to the field of hyperreals~$\astr$; see \cite{Ro66}.

In the 1970s, Hrbacek \cite{Hr78} and Nelson \cite{Ne77} independently
developed axiomatic approaches to NSA.  In such axiomatic approaches,
infinitesimals are found, not in an extension~$\astr$, but rather
in~$\R$ itself.  What enables this is the introduction of a one-place
predicate ``standard'' postulating the existence of two types of real
numbers (and other entities): standard and nonstandard.

In axiomatic approaches, an infinitesimal is a number smaller in
absolute value than every positive standard real number.  Hrbacek and
Nelson showed the conservativity of their theories over the
traditional Zermelo--Fraenkel set theory with the Axiom of Choice,
meaning that the introduction of the new predicate does not involve
any new foundational commitments and does not, in principle, prove any
new theorems.  If we were to go beyond our analysis of the Leibnizian
procedures, and ask which foundational approach to NSA would be closer
to the spirit of the Leibnizian calculus, we would probably have to
answer that the axiomatic approach is closer, since the existence of
inassignable infinitesimals is postulated by Leibniz, rather than
resulting from any kind of construction (as in Robinson's original
approach).

\end{document}